\documentclass[reqno]{amsart}
\usepackage{amsmath,amsfonts,amssymb,graphics,graphicx,latexsym,amscd,subfigure,hyperref}

\newtheorem{thm}{Theorem}[section]
\newtheorem{cor}[thm]{Corollary}
\newtheorem{prop}[thm]{Proposition}
\newtheorem{question}[thm]{Question}
\newtheorem{lemma}[thm]{Lemma}
\newtheorem{defn}[thm]{Definition}
\newtheorem{rem}[thm]{Remark}

\newtheorem{exa}[thm]{Example}

\newcommand{\bbQ}{\mathbb{Q}}
\newcommand{\bbR}{\mathbb{R}}
\newcommand{\bbT}{\mathbb{T}}
\newcommand{\bbC}{\mathbb{C}}
\newcommand{\bbZ}{\mathbb{Z}}
\newcommand{\ip}{\cdot}
\newcommand{\ismax}{\operatorname{Iso}^{\operatorname{max}}}
\newcommand{\is}{\operatorname{Iso}}

\newcommand{\g}{\gamma}

\newcommand{\opF}{\operatorname{Fix}}
\newcommand{\orb}{\mathcal{O}}
\newcommand{\order}{\Omega}

\newcommand{\tf}{\tilde{f}}
\newcommand{\tn}{\widetilde{N}}

\newcommand{\tx}{\widetilde{x}}
\newcommand{\ty}{\widetilde{y}}

\begin{document}

\title[]{Equivariant inverse spectral theory and toric orbifolds}

\author{Emily B. Dryden}
\address{Department of Mathematics, Bucknell University, Lewisburg, PA 17837, USA}
\email{ed012@bucknell.edu}
\author{Victor Guillemin}
\address{Department of Mathematics, Massachusetts Institute of Technology, Cambridge, MA 02139, USA}
\email{vwg@math.mit.edu}
\author{Rosa Sena-Dias}
\address{Centro de An\'alise Matem\'atica, Geometria e Sistemas Din\^amicos, Departamento de Matem\'{a}tica, Instituto Superior T\'{e}cnico, Av. Rovisco Pais, 1049-001 Lisboa, Portugal}
\email{rsenadias@math.ist.utl.pt}

\date{}

\begin{abstract}
Let $\orb^{2n}$ be a symplectic toric orbifold with a fixed $\mathbb{T}^n$-action and with a toric K\"ahler metric $g$.  In \cite{dgs} we explored whether, when $\orb$ is a manifold, the equivariant spectrum of   
the Laplace operator $\Delta_g$ on $\mathcal{C}^\infty(\orb)$ determines the moment polytope of $\orb$, and hence by Delzant's theorem determines $\orb$ up to symplectomorphism.  In the setting of toric orbifolds we significantly improve upon our previous results and show that the moment polytope of a generic toric orbifold is determined by its equivariant spectrum, up to two possibilities and up to translation.  This involves developing the asymptotic expansion of the heat trace on an orbifold in the presence of an isometry.  We also show that the equivariant spectrum determines whether the toric K\"{a}hler metric has constant scalar curvature.
\end{abstract}


\maketitle


\section{Introduction}\label{sec:intro}

Given a Riemannian manifold $(M, g)$, one can consider the Laplace operator $\Delta_g$ acting on the space of smooth functions on $M$; the spectrum of $\Delta_g$ is the set of eigenvalues of $\Delta_g$ on $\mathcal{C}^\infty(M)$.  From a spectral-theoretic point of view, one is interested in how much about the geometry of $(M,g)$ is determined by the spectrum of $\Delta_g$.  There are examples of Riemannian manifolds with the same spectrum which are not isometric (e.g., \cite{Gor_survey}, \cite{Mil_tori}, \cite{Sunada_method}), and there are also positive results showing that manifolds within a certain class are spectrally determined (e.g., \cite{Tanno}).  In the setting of symplectic toric geometry, Miguel Abreu \cite{m2} asked
\begin{question}\label{q:Abreu}
Let $M$ be a toric manifold equipped with a toric K\"ahler metric $g$. Does the spectrum of the Laplacian $\Delta_g$ determine the moment polytope of $M$?
\end{question}

In \cite{dgs} the authors considered a modified version of this question, replacing the spectrum of the 
Laplacian by the \emph{equivariant spectrum} of the Laplacian. This is simply the spectrum of the Laplacian together with, for each eigenvalue, the weights of the representation of $\bbT^n$ on the eigenspace corresponding to the given eigenvalue.  Question \ref{q:Abreu} then becomes
\begin{question}\label{q:ours}
 Let $M$ be a toric manifold equipped with a toric K\"{a}hler metric $g$. Does the equivariant spectrum of $\Delta_g$ on $\mathcal{C}^\infty(M)$ determine the moment polytope of $M$?
\end{question}

Thomas Delzant \cite{Delzant} proved that the moment polytope of a toric symplectic manifold $M$ determines $M$ up to symplectomorphim.  Thus, if the answer to Question \ref{q:ours} is ``yes,'' the spectrum of the Laplacian of a symplectic toric manifold determines its symplectomorphism type.  We showed that the answer is positive for many generic toric $4$-manifolds, up to translation and a small number of choices; however, we could not resolve the question when the moment polytope of $M$ has ``many'' parallel sides, a case which occurs with positive probability.

Toric orbifolds are a natural generalization of toric manifolds. They admit toric K\"{a}hler metrics, i.e., metrics that are determined by a symplectic form and a compatible, integrable almost complex structure and that are invariant under the torus action.  Thus one may again define the Laplacian and its associated equivariant spectrum. Toric orbifolds also have moment polytopes associated to them, so it is natural to ask Question \ref{q:ours} in the context of toric orbifolds.  The same issues as in \cite{dgs} arise when the moment polytope has parallel facets, but unlike for manifolds such facets occur for orbifolds with zero probability.  Hence we are able to prove the following theorem.

\begin{thm}\label{main_theorem}
Let $\orb$ be a generic toric orbifold with a fixed torus action and a toric K\"ahler metric. Then the equivariant spectrum of $\orb$ determines the moment polytope $P$ of $\orb$, and hence the equivariant symplectomorphism type of $\orb$, up to two choices and up to translation.
\end{thm}

\noindent Note that the two choices determined by the equivariant spectrum have symplectomorphic underlying manifolds (see \S \ref{sec:mainthm}).

The main tool in \cite{dgs} is a result of Harold Donnelly \cite{Don1} which gives an asymptotic expansion for the heat kernel in the presence of an isometry on a manifold. A key step in the proof of Theorem \ref{main_theorem} is to generalize this tool to the setting of orbifolds (see \S \ref{sec:eq_heat_oflds}). Our approach is similar to the generalization of the asymptotic expansion of the heat kernel to orbifolds that was done in \cite{dggw}; the resulting expansion should be of independent interest. 

\begin{thm}\label{main_asymptotic} 
Let $\orb$ be a closed Riemannian orbifold, let $K(t,x,y)$ be the heat kernel of $\orb$, and let $f$ be a nontrivial liftable isometry of $\orb$. Then $\int_\orb\,K(t,x,f(x))\text{dvol}_{\orb}(x)$ is asymptotic as $t\to 0^+$ to
\[
 \sum_{S \in \mathcal{S}(\orb)} \frac{1}{|\is(S)|} (4\pi t)^{-\frac{\text{dim}(\text{Fix} f \cap S)}{2}} \int_{\text{Fix} f \cap S} \sum_{k=0}^{\infty} b_k(f,S) t^k \text{dvol}_{\text{Fix} f \cap S}(x),
\]
where $S(\orb)$ is a stratification of $\orb$ and $|\is(S)|$ denotes the order of the isotropy of any point $p \in S$.  
\end{thm}

\noindent In a subsequent paper the authors will extend this result to more general operators than the Laplacian using semi-classical analysis techniques.

 Theorem \ref{main_theorem} holds for any toric K\"ahler metric on $\orb$, and it is well known (see \cite{g1}, \cite{m1}) that toric orbifolds admit many toric K\"ahler metrics. Thus one is led to ask what the equivariant spectrum tells us about the toric metric itself.
\begin{question}
Does the equivariant spectrum corresponding to a toric K\"ahler metric on a toric orbifold determine the toric K\"ahler metric?
\end{question}
A positive answer to this question is unlikely, but one could hope that the equivariant spectrum might determine some properties of the metric.  Finding ``special'' K\"ahler metrics on K\"ahler manifolds or orbifolds is currently an active research topic, with especial attention to K\"ahler Einstein metrics and extremal metrics in the sense of Calabi.  As a particular instance of these, one often looks for constant scalar curvature metrics. It is known that such metrics do not always exist, but it is conjectured that their existence is equivalent to a stability condition on the underlying manifold. For $4$-dimensional toric manifolds this conjecture was proved recently by Simon Donaldson \cite{Doncsc}. We will use the asymptotic expansion in Theorem \ref{main_asymptotic} to show that one can equivariantly hear constant scalar curvature toric K\"ahler metrics. 
\begin{thm}\label{thm:csc}
Let $\orb$ be a generic toric orbifold endowed with a toric K\"ahler metric $g$. Then the equivariant spectrum of $\orb$ determines if $g$ has constant scalar curvature.
\end{thm}
 
The paper is organized as follows.  In \S \ref{sec:ofld_back} we give the necessary background on orbifolds, with particular emphasis on orbifold strata and isotropy groups.  This allows us to prove Theorem \ref{main_asymptotic} in \S \ref{sec:eq_heat_oflds}.  We then specialize to the setting of toric orbifolds, giving relevant background in \S \ref{sec:toric_back}.  The proof of Theorem \ref{main_theorem} is given in \S \ref{sec:mainthm}, followed by the proof of Theorem \ref{thm:csc} in \S \ref{sec:csc}.
 
 \vspace{.5cm}
\noindent \textbf{Acknowledgments:}  The first and third authors appreciate the hospitality shown to them by the Mathematics Department at MIT during their visits there.  The first author's visit during Summer 2010 was partially supported by an NSF-AWM Mentoring Travel Grant, and the second author was partially supported by NSF grant DMS-1005696. We thank Yael Karshon and Isabella Novik for making us aware of the work by Daniel Klain on the Minkowski problem.
 

\section{Background on orbifolds}\label{sec:ofld_back}

We begin by reviewing some of the basic definitions related to orbifolds that are relevant to our work.  Our presentation and notation will follow that used in \S 2 of \cite{dggw}, which the reader may consult for more details.

A $k$-dimensional \emph{orbifold} $\orb$ is a second-countable Hausdorff topological space $X$ that is equipped with a maximal orbifold atlas.  Each chart in this atlas consists of a connected open subset $\tilde{U} \subset \mathbb{R}^k$, a finite group $G_U$ acting on $\tilde{U}$ by diffeomorphisms, and a mapping $\pi_U: \tilde{U} \rightarrow U$, where $U$ is an open subset of $X$ and $\pi_U$ induces a homeomorphism from $G_U \backslash \tilde{U}$ onto $U$.  We will assume that the action of $G_U$ on $\tilde{U}$ is effective.

Points in $\orb$ are either \emph{singular} or \emph{regular}.  A point $x \in \orb$ is singular if for some (hence every) orbifold chart $(\tilde{U}, G_U, \pi_U)$ about $x$, the points in the inverse image of $x$ in $\tilde{U}$ have nontrivial isotropy in $G_U$.  The isomorphism class of this isotropy group is called the \emph{abstract isotropy type} of $x$ and is independent of the choice of point in the inverse image of $\pi_U$ and of the choice of orbifold chart about $x$.    

An orbifold $\orb$ can be endowed with a \emph{Riemannian structure} by assigning to each orbifold chart $(\tilde{U}, G_U, \pi_U)$ a $G_U$-invariant Riemannian metric on $\tilde{U}$ satisfying an appropriate compatibility condition among charts.  Every Riemannian orbifold has an associated orthonormal frame bundle, as we now briefly describe.  For an orbifold $\orb$ of the form $G \setminus M$, where $M$ is a Riemannian manifold and $G$ is a discrete subgroup of the isometry group of $M$,  we begin by considering the orthonormal frame bundle of $M$, $F(M) \rightarrow M$.  Since each element of $G$ induces a diffeomorphism of $F(M)$ that takes fibers to fibers, we have an action of $G$ on $F(M)$ that covers the action of $G$ on $M$.  We then define the orthonormal frame bundle of $O$, $F(\orb)$, to be $G \setminus F(M) \rightarrow \orb$, with the fiber over a point $x \in \orb$ defined as the preimage of $x$ in $G \setminus F(M)$.  Note that $F(M)$ admits a right action of the orthogonal group $O(k)$, and this action commutes with the left action of $G$; hence we have a right $O(k)$-action on $F(\orb)$.  For an arbitrary orbifold $\orb$, one may define an orthonormal frame bundle using the above construction on the local charts on $\orb$ with an appropriate compatibility condition among charts; this bundle is an \emph{orbibundle} whose total space is again a smooth manifold that admits a right action of the orthogonal group.  The orbifold $\orb$ can then be viewed as the orbit space $F(\orb) / O(k)$.      

Let $x \in \orb$ be a singular point and view it as an element of $F(\orb)/ O(k)$.  As $\tilde{x} \in F(\orb)$ ranges over the fibers in the preimage of $x$ in $F(\orb)$, the stabilizer $\is (\tilde{x})$ ranges over a conjugacy class of subgroups of $O(k)$.  It turns out that this conjugacy class is independent of the choice of Riemannian metric used to construct $F(\orb)$, so it makes sense to call it the \emph{isotropy type} of $x \in \orb$, denoted $\is(x)$.  Moreover, the subgroups in a given conjugacy class can be shown to lie in the isomorphism class defined by the abstract isotropy type of $x$.  The cardinality $|\is(x)|$ is the \emph{order} of the isotropy at $x$, and the equivalence classes of points with the same isotropy type are called the \emph{isotropy equivalence classes}. 

\begin{defn}
A smooth \emph{stratification} of an orbifold $\orb$ is a locally finite partition of $\orb$ into  submanifolds.  Each submanifold, called a \emph{stratum}, is locally closed and its closure is the union of the stratum with a collection of lower-dimensional strata.
\end{defn}

Given a stratification of an orbifold $\orb$, the strata of maximal dimension are open in $\orb$ and their union has full measure in $\orb$.  From general results about smooth actions of Lie groups on manifolds, one can prove

\begin{prop}\cite[Cor. 2.11]{dggw}
Let $\orb$ be an orbifold.  Then the action of $O(k)$ on the frame bundle $F(\orb)$ gives rise to a (Whitney) stratification of $\orb$.  The strata are connected components of the isotropy equivalence classes in $\orb$.  The set of regular points of $\orb$ intersects each connected component $\orb_0$ of $\orb$ in a single stratum that constitutes an open dense submanifold of $\orb_0$.
\end{prop}

The strata of $\orb$ will be called \emph{$\orb$-strata}; given an orbifold chart $(\tilde{U}, G_U, \pi_U)$ on $\orb$, the stratifications of sets $U$ and $\tilde{U}$ induced by $G_U$ will be referred to as \emph{$U$-strata} and \emph{$\tilde{U}$-strata}, respectively.  The following result was proved in \cite{dggw}.

\begin{prop}\cite[Prop. 2.13]{dggw}
Let $\orb$ be a Riemannian orbifold and
$(\tilde{U},G_U,\pi_U)$  be an orbifold chart.  Then:
 \begin{enumerate}
 \item The $U$-strata are precisely the connected components of the intersections of
the $\orb$-strata with $U$.  
\item Any two elements of the same $\tilde{U}$-stratum have the same stabilizers
 in $G_U$ (not just conjugate stabilizers).  
\item  If $H$ is a subgroup of $G_U$, then each connected component $W$ of the
fixed point set $\opF(H)$ of $H$ in $\tilde{U}$ is a closed submanifold of $\tilde{U}$.  Any
$\tilde{U}$-stratum  that intersects $W$  nontrivially lies entirely in $W$.  Thus the
stratification of $\tilde{U}$ restricts to a stratification of $W$.
\end{enumerate}
\end{prop}
It follows from this proposition that if $\tilde{N}$ is any $\tilde{U}$-stratum in $\tilde{U}$, then all the points in $\tilde{N}$ have the same isotropy group in $G_U$; we denote this isotropy group of $\tilde{N}$ by $\is(\tilde{N})$.  The set of $\gamma \in \is(\tilde{N})$ such that $\tilde{N}$ is open in $\opF(\gamma)$ will be denoted $\is^{\operatorname{max}} (\tilde{N})$.  Note that the union of the $\tilde{U}$-strata $\tilde{N}$ for which $\gamma$ is an element of $\is^{\operatorname{max}} (\tilde{N})$ has full measure in $\opF(\gamma)$. 

\begin{exa}\label{exa1}
Let $M \subset \mathbb{R}^3$ be a sphere of radius $1$ centered at the origin.  The quotient of $M$ by a rotation of order $p$ about the $z$-axis is a closed orbifold called a $(p,p)$-football.  Each of $(1,0,0)$ and $(0,0,-1)$ is a singular stratum with isotropy of order $p$. 

We may also consider a $(p,q)$-football with $p \neq q$; this does not arise as a global quotient of a manifold, but is an orbifold with two antipodal singular points whose underlying topological space is a sphere.  Taking the quotient of this orbifold with respect to reflection in the plane containing the origin and the singular points, we get an orbifold whose underlying space is a disk.  The points on the boundary of the disk are singular points (not boundary points) of the orbifold and comprise four strata: the image of each antipodal point forms a single stratum with isotropy of order $2p$ or $2q$, respectively,  and each open edge forms a stratum with isotropy of order $2$.  The intersection $U$ of the disk with a smaller disk centered at one of the ``poles'' is the image of an orbifold chart $(\tilde{U}, G, \pi_U)$, where $\tilde{U}$ is a disk in $\mathbb{R}^2$ centered at the origin and $G$ is the dihedral group of appropriate order.  The $\tilde{U}$-strata with respect to $G$ are the origin and the intersection of $\tilde{U}$ with, say, the positive and negative $x$-axis.  Let $\tilde{N}$ be the intersection of $\tilde{U}$ with one of the half-axes; then $\is(\tilde{N})$ contains the reflection and the identity, but $\ismax(\tilde{N})$ is just the reflection.  If $\tilde{N} = \{0\}$, then $\is(\tilde{N})=G$ but $\ismax(\tilde{N})$ is just the rotation through angle $\frac{2 \pi}{p}$ (respectively, $\frac{2\pi}{q}$) about the origin. 
\end{exa}
\noindent We will return to variations on this example in \S \ref{sec:toric_back}.


\section{Equivariant heat kernel asymptotics for orbifolds}\label{sec:eq_heat_oflds}

We now develop the asymptotic expansion of the heat kernel on an orbifold $\orb$ in the presence of a liftable isometry.  To show that such an expansion exists and to find it, one approach is to use the local structure of orbifolds.  In particular, one may take a local covering of $\orb$ by convex geodesic balls and piece together a parametrix for the equivariant heat operator on $\orb$ using locally defined parametrices: working in a convex geodesic ball $U \subset \orb$ with orbifold chart $(\tilde{U}, G_U, \pi_U)$, define a local parametrix $\tilde{H} (t,\tilde{x},\tilde{y})$ on $(0, \infty) \times \tilde{U} \times \tilde{U}$.  Suppose that $f: \orb \rightarrow \orb$ is an isometry of $\orb$ that lifts to an isometry $\tilde{f}: \tilde{U} \rightarrow \tilde{U}$ with $\tilde{f} \circ \gamma = \gamma \circ \tilde{f}$ for all $\gamma \in G_U$.  Then the function
\[
(t, \tilde{x}, \tilde{y}) \longmapsto \sum_{\gamma \in G_U} \tilde{H} (t, \tilde{x}, \tilde{f} \circ \gamma (\tilde{y}))
\]
descends to a well-defined function on $(0, \infty) \times U \times U$.  The argument to patch together these locally defined parametrices to get a globally defined parametrix and thus an equivariant heat kernel on $\orb$ follows as in \S 3 of \cite{dggw}.

To find the asymptotic expansion of the equivariant heat kernel, we generalize the following theorem of Donnelly.
\begin{thm}\label{prop:Don}\cite{Don1}.  Let $M$ be a closed Riemannian manifold,
let $K(t,x,y)$ be the heat kernel of $M$, and let $\g$ be a nontrivial isometry
of $M$. 
Then 
$\int_M\,K(t,x,\g(x))\text{dvol}_M(x)$ is asymptotic as $t\to 0^+$ to
$$\sum_{W\subset\opF(\g)}\,(4\pi
t)^{-\frac{\dim(W)}{2}}\sum_{k=0}^\infty\,t^k\int_W\,b_k(\g,a)\text{dvol}_W(a)$$
where $W$ ranges over connected components of the fixed point set of $\gamma$, $b_k(\gamma,a)$ is a real-valued function on the fixed point set of $\g$
and $\text{dvol}_W$ is the volume form on $W$ defined by the Riemannian metric
induced from $M$.
 \end{thm}

The function $b_k(\g, x)$ has several key properties (see \cite[\S 4]{dggw}).  First, its restriction to any $W  \subset \opF(\g)$ is smooth.  Second, it is local in that it only depends on the germs at $x$ of the Riemannian metric and of the isometry $\g$.  Finally, it is universal in that it behaves as one would hope with respect to isometries; namely, if $M$ and $M'$ are Riemannian manifolds admitting the isometries $\g$ and $\g'$, respectively, and $\sigma: M \rightarrow M'$ is an isometry satisfying $\sigma \circ \g = \g' \circ \sigma$, then $b_k(\g,x) = b_k (\g', \sigma(x))$ for all $x \in \opF (\g)$.  Donnelly gave explicit formulas for $b_0$ and $b_1$, with the general definition of $b_k$ as follows.  Let $x \in W$.  Note that the orthogonal complement of $T_xW$ in $T_xM$ is invariant under $\gamma_*$; let $A_\g(x)$ be the nonsingular matrix transformation defined by $\g_*$, and set $B_\g(x) = (I - A_\g(x))^{-1}$.  Then we can define
\[
b_k(\g, x) = |\det (B_\g(x))| b'_k(\g, x),
\]
where $b'_k (\g, x)$ is a universal invariant polynomial in the components of $B_\g$ and in the curvature tensor $R$ of $M$ and its covariant derivatives at $x$. 

With this definition in mind, we prove a special case of Theorem \ref{main_asymptotic}.

\begin{lemma}\label{lemma:global}
Let $\orb=G\backslash M$, where $M$ is a Riemannian manifold and $G$ is a finite group acting effectively on $M$.  Let $f: \orb \rightarrow \orb$ be an isometry that lifts to an isometry $\tilde{f}: M \rightarrow M$ with $\tilde{f} \circ \gamma = \gamma \circ \tilde{f}$ for all $\gamma \in G$.  Then 
\[
\sum_{i=1}^{\infty}\operatorname{Tr} (f_{\lambda_i})^* e^{-t\lambda_i} \sim 
\sum_{S \in \mathcal{S}(\orb)} \frac{1}{|\is(S)|} (4\pi t)^{-\frac{\dim(\opF f \cap S)}{2}} \sum_{k=0}^{\infty} t^k \int_{\opF f \cap S} b_k(f,x) \text{dvol}_{\opF f \cap S}(x)
\]
as $t \to 0^+$, where $|\is(S)|$ is the order of the isotropy at every point in $S$ as defined in \S \ref{sec:ofld_back}.
\end{lemma}

\begin{proof}
Note that $(M,G,\pi)$ is a global orbifold chart where
$\pi:M\to\orb$ is the projection.  If $K$ denotes the heat kernel of $M$, then the heat kernel $K^\orb$ of $\orb$ is given by
$$K^\orb(t,x,y)=\sum_{\g\in G}\,K(t,\tx,\g(\ty))$$
where $\tx$, respectively $\ty$, are any elements of $\pi^{-1}(x)$, respectively
$\pi^{-1}(y)$.  Thus 
\[
\int_\orb K^\orb(t,x,x) \text{dvol}_\orb(x)=\frac{1}{|G|}\sum_{\g\in
G}\,\int_M\,K(t,\tx,\g(\tx))\text{dvol}_M(\tx).
\] 

Let us examine what happens to this heat kernel in the presence of an isometry $f: \orb \rightarrow \orb$.  Letting $K_f^{\orb}$ denote the $f$-equivariant heat kernel of $\orb$, we have that 
\[
\int_\orb K_f^\orb(t,x,x) \text{dvol}_\orb(x)=\frac{1}{|G|}\sum_{\g\in
G}\,\int_M\,K(t,\tx,\tf (\g(\tx)))\text{dvol}_M(\tx).
\]
Applying Theorem \ref{prop:Don} to this expression gives
\begin{equation}\label{eqn:Dondirect}
\int_\orb K_f^\orb(t,x,x) \text{dvol}_\orb(x) \sim \frac{1}{|G|}\sum_{\g\in G} \sum_{W \in \opF(\tf \circ \gamma)} (4\pi
t)^{-\frac{\dim(W)}{2}}\sum_{k=0}^\infty\,t^k\int_W\,b_k(\tf \circ \g,\tx)\text{dvol}_W(\tx)
\end{equation}
as $t\to 0^+$.
In order to express the right side of \eqref{eqn:Dondirect} in terms of intrinsic orbifold data, we need to analyze the fixed point sets arising from $\tf \circ \gamma$ for $\gamma \in G$.  We begin by relating these fixed point sets to $M$-strata. 

Let $W$ be a connected component in $\opF(\tf \circ \g)$ and let $\tn$ be an $M$-stratum contained in $W$.  Then either $\tn$ has measure zero in $W$ (in which case
$\g\notin\ismax(\tn)$) or $\tn$ is open in $W$ and $\g \in \ismax(\tn)$.  
Suppose $\tx \in \tn$ and $\g \in \ismax(\tn)$.  Since $\tn \subset W \in \opF(\tf \circ \g)$, we have $\tf(\g(\tx)) = \tx$.  But $\g \in \ismax(\tn)$, so $\g(\tx)=\tx$.  Thus $\tf(\tx) = \tx$, or $\tx \in \opF(\tf \circ I)$. This means that we only need to consider the contribution from the identity element in $G$, and we may replace the integral over $W$ with integrals over the $M$-strata that are open in $W$.  Thus the right side of \eqref{eqn:Dondirect} becomes
\begin{equation}\label{eqn:trace_up}
\frac{1}{|G|} \sum_{\tn \in \mathcal{S}_f(M)} (4 \pi t)^{-\frac{\dim(\tn)}{2}} \sum_{k=0}^{\infty} t^k \int_{\tn} b_k(\tf \circ \g,\tx)\text{dvol}_{\tn}(\tx),
\end{equation}
where $\mathcal{S}_f(M)$ denotes the strata in $\mathcal{S}(M)$ that are open in $\opF(\tf \circ \g)$. 
   
Our next task is to relate the data on the manifold $M$ to data on $\orb$.  Let $N$ be an $\orb$-stratum that is open in a component of $\opF(f)$.  Then $\pi^{-1}(N)$ is a union of finitely many mutually isometric strata in $\mathcal{S}_f(M)$ and $\pi: \pi^{-1}(N) \to N$ is a covering map of degree $\frac{|G|}{|\operatorname{Iso}(N)|}$.  Moreover, the total contributions to \eqref{eqn:trace_up} from the elements of $\pi^{-1}(N)$ are equal to 
\[
\frac{|G|}{|\operatorname{Iso}(N)|} (4 \pi t)^{-\frac{\dim(N)}{2}} \sum_{k=0}^{\infty} t^k \int_{N} b_k(f, x)\text{dvol}_{N}(x).
\]
Thus \eqref{eqn:trace_up} becomes
\begin{equation}\label{eqn:tracedown}
 \sum_{N \in \mathcal{S}_f(\orb)}\frac{1}{|\operatorname{Iso}(N)|} (4 \pi t)^{-\frac{\dim(N)}{2}} \sum_{k=0}^{\infty} t^k \int_{N} b_k(f, x)\text{dvol}_{N}(x)
\end{equation}
where $\mathcal{S}_f(\orb)$ denotes the strata in $\mathcal{S}(\orb)$ that are open in $\opF(f)$.  This proves the lemma.
\end{proof}

The argument in the proof of Lemma \ref{lemma:global} can be applied to orbifold charts, and one may piece together the resulting computations via a partition of unity to prove Theorem \ref{main_asymptotic}.  The ideas are exactly the same as those used in \cite{dggw} to pass from the asymptotic expansion of the heat trace for an orbifold of the form $\orb=G\backslash M$ to the expansion for a general orbifold.  We refer the interested reader to \S 4 of \cite{dggw} for details.

Using computations from \cite{Don1} and \cite{dggw} one can find explicit expressions for the first few terms in the asymptotic expansion in Theorem \ref{main_asymptotic}. We will denote the scalar curvature by $s$, the Ricci tensor by $\rho$ and the full curvature tensor by $R$. Let $S$ be a connected component of $\mathcal{S}(\orb)$ and $x\in S$. Let $\gamma$ be an element of $\ismax(S)$ and $W_f$ be a local lift of $S\cap \opF(f)$ via an orbifold chart. Let $A_{f,\gamma}$ be the isometry 
\[
 df\circ d\gamma :{T_{\tx}W_f^\perp}\rightarrow T_{\tx}W_f^\perp,
\]
where $T_{\tx}W_f^\perp$ denotes the normal space to $W_f$, and set $B_{f,\gamma}(\tx)=(I-A_{f,\gamma}(\tx))^{-1}$. Then we have
\[
 b_0(f,x)=\sum_{\gamma \in \ismax(S)}|\det (B_{f,
\gamma}(\tx))|,
\]
and
\[
 b_1(f,x)=\sum_{\gamma \in \ismax(S)}|\det (B_{f,\gamma}(\tx))|\tau_\gamma(\tx)
\]
where 
\[
 \tau_\gamma=\frac{s}{6}+\frac{1}{6}\rho_{kk}+\frac{1}{3}R_{iksh}
B_{ki}B_{hs}+\frac{1}{3}R_{ikth}
B_{kt}B_{hi}-R_{kaha}
B_{ks}B_{hs}.
\]
Here the indices $k,i,s,h,t$ correspond to normal directions in $T_{\tx}W_f^\perp$ and we sum over repeated indices. 

When $f$ is the identity and we are in a neighborhood of a regular point $x \in \orb$, the function $b_k$ equals the usual heat invariant $a_k$ for manifolds.  In particular, we note for use in \S \ref{sec:csc} that
\begin{equation}\label{eqn:b2}
 b_2(x)=\frac{1}{360} \left(2|R|^2-2|\rho|^2+5s^2 \right).
\end{equation}


\section{Background on toric orbifolds} \label{sec:toric_back}

We now specialize to the setting of toric orbifolds, providing the definitions and background that are needed to understand the proof of Theorem \ref{main_theorem}.
The notion of symplectic manifold generalizes easily to the orbifold setting: an orbifold is said to be symplectic if it admits a $2$-form $\omega$ which is nondegenerate and closed.  One can then specify what it means for an orbifold to be \emph{toric}.

\begin{defn}
Let $(\orb,\omega)$ be a symplectic orbifold of real dimension $2n$. Then $(\orb,\omega)$ is said to be toric if admits an effective Hamiltonian $\bbT^n$-action, where $\bbT^n$ is the real torus of dimension $n$.
\end{defn}
An action of a Lie group $G$ on a symplectic orbifold $(\orb,\omega)$ is said to be \emph{Hamiltonian} if it admits a \emph{moment map}.  This is a map $\phi:\orb \rightarrow \mathfrak{G}^*$, where $\mathfrak{G}^*$ denotes the dual of the Lie algebra of $G$, satisfying
\[
d\phi(x)(v)=v^\sharp \lrcorner \omega,
\]
where $v$ is an element of the Lie algebra of $G$ and $v^\sharp$ is the vector field determined by $v$ on $\orb$.  That is, if $g(t)$ is a smooth path on $G$ with $g(0)=\text{id}$ and
\[
v=\frac{d}{dt}_{|t=0}g(t),
\]
then
\[
v^\sharp(x)=\frac{d}{dt}_{|t=0}g(t)\cdot x.
\]
Note that this is well-defined up to a constant.

For toric manifolds it is known that the image of the moment map determines the manifold up to symplectomorphism. For toric orbifolds the image of the moment map is insufficient to give this result, but Eugene Lerman and Sue Tolman \cite{LT} have indicated what additional data one needs to determine the orbifold.  To state their result, we begin by describing the image of the moment map.
Identifying the dual of the Lie algebra of $\bbT^n$ with $\bbR^n$, Lerman and Tolman showed that the image of the moment map of a toric orbifold is a special type of convex polytope in $\bbR^n$.  
\begin{defn}\label{defn:rat_simple}
 A convex polytope $P$ in $\mathbb{R}^n$ is \emph{rational simple} if
 \begin{enumerate}
 \item there are $n$ facets meeting at each vertex;
 \item for every facet of $P$, a primitive outward normal can be chosen in $\bbZ^n$;
  \item for every vertex of $P$, the outward normals corresponding to the facets meeting at that vertex form a basis for $\bbQ^n$.
\end{enumerate}
Note that a facet is a face in $P$ of codimension $1$.
\end{defn}
\noindent The triangle with vertices $(0,0), (0,1), (2,0)$  is an example of a rational simple polytope.  Moreover, this polytope cannot be the moment map image of a toric manifold; in the manifold setting, we replace $\bbQ^n$ in (3) of Definition \ref{defn:rat_simple} by $\bbZ^n$.  

It is also important to understand what types of singularities can occur in toric orbifolds.  Lerman and Tolman show that all the points over the interior of the moment polytope of the toric orbifold are regular and that the points over the facets of the polytope have cyclic isotropy type. 
\begin{thm}\cite[Theorem 6.4]{LT}\label{thm:isotropy}
Let $(\orb,\omega)$ be a toric orbifold with moment map $\phi$. Let $F$ be an open facet of $\phi(\orb)$. Then there exists an integer $m_F$ such that all points in $\phi^{-1}(F)$ have isotropy group $\bbZ_{m_F}$. 
\end{thm}
The integer $m_F$ is called the label of the open facet $F$, and we call the moment polytope together with the facet labels a \emph{labeled polytope}. Thus to each toric orbifold of dimension $2n$, we can associate a labeled rational simple polytope in $\bbR^n$.  Lerman and Tolman proved that these labeled polytopes essentially determine the associated symplectic orbifolds. 
\begin{prop}\cite[Proposition 6.5]{LT}
If two toric orbifolds have the same labeled moment polytopes up to $SL(n,\bbZ)$-transformations and translations, then the orbifolds are equivariantly symplectomorphic.
\end{prop}

\begin{exa}
We return to the $(p,q)$-footballs of Example \ref{exa1}.  These are symplectic orbifolds and they admit a Hamiltonian $S^1$-action given by rotation about the north-south axis.  The labeled polytope associated to this toric orbifold is the interval $[-1,1]$ with labels $p$ and $q$ at the upper and lower endpoints, respectively.
\end{exa}

In the proof of Theorem \ref{main_theorem}, we will need two further results relating the moment polytope to its associated orbifold.  These results are stated and proved in the setting of toric manifolds in \cite{dgs}, but the proofs are exactly the same in the orbifold case.  First we examine the fixed point set of an isometry of an orbifold.
\begin{lemma}\cite[Lemma 2.9]{dgs}\label{fixed_point_set}
 Let $\theta\in \bbR^n$. The fixed point set of $\psi(e^{i\theta})$, denoted $F_\theta$, is the union of the pre-images via the moment map of all faces to which $\theta$ is normal in a face of lower codimension.
\end{lemma}

Finally, we give the relationship between the volume of a face in the polytope and the volume of its pre-image under the moment map.

\begin{lemma}\cite[Lemma 2.10]{dgs}\label{vol}
Consider a face $F$ of dimension $q$ in the labeled polytope $P$ of a symplectic toric orbifold $(\orb, \omega)$.  Let $\phi$ be the moment map of the torus action with respect to the form $\omega$.  Then
\[
\operatorname{Vol}_{\omega} (\phi^{-1}(F)) = (2 \pi)^q \operatorname{Vol}(F).
\]
\end{lemma}

Since we want to study the spectrum of the Laplacian on toric orbifolds, we need to understand metrics on such orbifolds. We will restrict to metrics that are compatible with the symplectic structure and that are K\"ahler.  That is, we consider metrics that come from an integrable almost complex structure $J$ on $(\orb,\omega)$ which is compatible with $\omega$; more precisely,
\[
g(\cdot,\cdot)=\omega(\cdot,J\cdot)
\]
defines a positive definite metric on $\orb$.

It was shown in \cite{g1} that all toric manifolds admit a special K\"ahler structure which is invariant under the torus action, called the reduced K\"ahler structure. We will refer to K\"ahler structures which are invariant under the torus action as toric K\"ahler structures. In \cite{m2} Abreu showed how to construct all other toric K\"ahler structures from the reduced K\"ahler structure using functions on the moment polytope of the toric manifold.  One may generalize the results in \cite{g1} and \cite{m2} to orbifolds; for a discussion of this generalization, see \cite{m3}.
\begin{thm}\cite{g1},\cite{m2}
Any toric orbifold admits many toric K\"ahler structures.
\end{thm}
As mentioned in \S \ref{sec:intro} the problem of finding the ``best" such K\"ahler structures has been the source of much work in differential geometry, and the K\"ahler structures which correspond to metrics whose scalar curvature is constant are of particular interest.  We return to these metrics in \S \ref{sec:csc}.


\section{Hearing a generic toric orbifold}\label{sec:mainthm}

In this section we prove Theorem \ref{main_theorem}. First we give a precise definition of equivariant spectrum; it is entirely analogous to the corresponding definition for manifolds.
\begin{defn}
Let $\orb^{2n}$ be a toric orbifold with a fixed torus action. Denote by $\psi:\bbT^n\rightarrow Sympl(\orb)$ the corresponding group homomorphism, and let $g$ be a toric metric on $\orb$. The \emph{equivariant spectrum} is the list of all the eigenvalues of the Laplacian on $(\orb,g)$  together with the weights of the action induced by $\psi(e^{i\theta})$ on the corresponding eigenspaces, for all $\theta \in \bbR^n$.  The eigenvalues and weights are listed with multiplicities. 
\end{defn}

By studying the asymptotic expansion in Theorem \ref{main_asymptotic}, we see that the equivariant spectrum provides significant information about the moment polytope of a toric orbifold.  We begin with the case in which the moment polytope has no parallel facets.
\begin{prop}\label{spectral_data}
 The equivariant spectrum associated to a toric orbifold $\orb$ whose moment polytope has no parallel facets determines
\begin{enumerate}
 \item the (unsigned) normal directions to the facets;
 \item the volumes of the corresponding facets;
\item the labels of the facets.
\end{enumerate} 
\end{prop}
\begin{proof}
Let $\phi$ be the moment map of the torus action on $\orb$ and let $P$ be its image, with $P$ given by
\begin{displaymath}
 P=\{x\in \bbR^n: x\ip u_i\geq c_i,\,\, i=1,\ldots,d\}
\end{displaymath}
for some collection of $u_i$ in $\bbR^n$ and $c_i \in \bbR$.
For each $u$ in $\bbR^n$ let $\psi(u)$ denote the isometry of $\orb$ given by the $e^{u}$-action on $\orb$. For each $\lambda$ in the spectrum of $\orb$, the map $\psi(u)$ induces a linear action on the $\lambda$-eigenspace of $\orb$ which we denote by $\psi^\sharp_\lambda(u)$.  The asymptotic expansion from Theorem \ref{main_asymptotic}  gives
\begin{equation}\label{eqn:asymp}
\sum_\lambda \text{tr} ( \psi^\sharp_\lambda(u))e^{-t\lambda}  \simeq 
\sum_{V_i} \frac{1}{|\is(V_i)|}(4\pi t)^{-\frac{\dim (V_i)}{2}}\text{Vol}(V_i)D_i +O(t) \ ,
\end{equation}
where the $V_i$ are the connected components of the fixed point set of $\psi(u)$ and $D_i$ is calculated as follows. Let $x$ be in $V_i$ and choose an orbifold chart $(\tilde{U}, G_U, \pi_U)$ about $x$, where $G_U = \is(x)$.  We can locally lift $\psi(u)$ to an isometry $\tilde\psi(\tilde{u})$ on $\tilde{U}$ which commutes with $G_U$. Then
\[
D_i(x)=\sum_{\g\in \ismax(V_i)} \left| \det((I-A_{\psi, \gamma}(\tilde{x}))^{-1}) \right|
\]
where $A_{\psi, \gamma}(\tilde{x})$ is as in \S \ref{sec:eq_heat_oflds}.  Note that $D_i$ is $G_U$-invariant and locally constant, so it is indeed a constant function on $V_i$. 
 
We will now untangle the polytope data contained in the right side of \eqref{eqn:asymp}. 
Set $u=ru_i$ for some $r$ and some $u_i$. If $r\ne 0$, Lemma \ref{fixed_point_set} implies that the connected components of the fixed point set of $\psi(u)$ are the pre-image via the moment map of all the faces whose normal is parallel to $u_i$. The highest-dimensional connected components have dimension $2(n-1)$ and we can determine the $u_i$'s up to sign: they are the vectors for which the right side of \eqref{eqn:asymp} has power $-(n-1)$ in $t$. When our polytope does not have parallel facets there is a single highest-dimensional connected component of $\text{Fix}(\psi(u))$, namely the pre-image of a facet $F$; we can determine its volume and its label by considering how $D_r$ varies with $r$. More precisely, given $x \in \phi^{-1}(F)$ we can find an orbifold chart $(\tilde{U}, G_U, \pi_U)$ about $x$ and we know by Theorem \ref{thm:isotropy} that $G_U$ is a cyclic group, say of order $\order$.  Let $\tilde{\psi}_r$ be the local lift of $\psi(ru_i)$. We have 
\[
D_r(x)=\sum_{l=1}^{\order-1} \left|\det((I-A_{\tilde{\psi_r}, \gamma^l}(\tilde{x}))^{-1})\right| 
\]
where $\gamma$ generates $G_U$. Now $T_{\tilde{x}}\operatorname{Fix}(\tilde{\psi_r})^{\perp}$ is a two-dimensional vector space and $A_{\tilde{\psi_r}, \gamma^l}$ are isometries of that space, so they must be rotations. Note that $d\gamma^{\order}$ is the identity and therefore the rotation angle of $d\gamma$ is a multiple of $2\pi/\order$; by reordering we may assume that the rotation angle equals $2\pi/\order$. 
Since $\tilde\psi_r$ is a group homomorphism from $S^1$ to the group of isometries of $\orb$, we see that $d\tilde\psi_r$ is a rotation whose angle is $r\theta$ for some fixed $\theta$. Therefore $A_{\tilde{\psi_r}, \gamma^l}$ is a rotation of angle $r\theta+\frac{2l\pi}{\order}$, and the matrix representation of $I-A_{\tilde{\psi_r}, \gamma^l}$ is 
\[
\begin{pmatrix}
1-\cos(r\theta+\frac{2l\pi}{\order})&-\sin(r\theta+\frac{2l\pi}{\order})\\
\sin(r\theta+\frac{2l\pi}{\order})&1-\cos(r\theta+\frac{2l\pi}{\order})
\end{pmatrix}.
\]
Thus 
\[
D_r=\sum_{l=1}^{\order-1} \frac{1}{2-2\cos(r\theta+\frac{2l\pi}{\order})},
\]
and the coefficient corresponding to the lowest-order term in the right side of \eqref{eqn:asymp} is  
\[
\frac{(4 \pi)^{-(n-1)}}{\order}\text{Vol}(F)\sum_{l=1}^{\order-1} \frac{1}{2-2\cos(r\theta+\frac{2l\pi}{\order})}.
\]
Hence this coefficient is spectrally determined for all $r\ne 0$, which implies that $\text{Vol}(F)$ and $\order$ are spectrally determined. Note that we have used that the volume of $\phi^{-1}(F)$ is proportional to the volume of $F$, as indicated in Lemma \ref{vol}.
\end{proof}

We next address the fact that the spectrum only determines the normals of the facets up to sign.  Note that a convex polytope with associated facet normals and volumes  $\{(u_i,\nu_i),i=1,\dots,d\}$ always satisfies 
\begin{displaymath}
\sum_{i=1}^d \nu_i u_i=0.
\end{displaymath}
\begin{defn}
Let $P$ be a convex polytope in $\bbR^n$ with associated facet normals and volumes  $\{(u_i,\nu_i),i=1,\ldots,d\}$. We say that $P$ has no subpolytopes if
\begin{displaymath}
 \sum_{i\in I}\nu_i u_i=0 \ \ \text{implies} \ \ I=\{1,\ldots,d \} \,\,\text{or}\,\, I=\emptyset.
\end{displaymath}
\end{defn}

For convex polytopes with no subpolytopes the set of normals up to sign determines the actual normals up to a finite number of sign choices. 
\begin{lemma}\label{two_choices}
Let $P$ be a convex polytope in $\bbR^n$ with no subpolytopes and facet volumes $\nu_1,\dots,\nu_d$. Assume that the facet normals to $P$ are $u_1, \dots, u_d$ up to sign. Then, up to translation, there are only $2$ choices for the set of signed normals.
\end{lemma}
\begin{proof}
For each choice $\xi_i=\pm u_i, i=1,\dots , d$, corresponding to a convex polytope we have
\begin{displaymath}
\sum_{i=1}^d \nu_i \xi_i=0.
\end{displaymath}
Adding the sums corresponding to two different choices, we get a relation 
\begin{displaymath}
\sum_I \nu_i \xi_i=0.
\end{displaymath}
Our assumption that $P$ has no subpolytopes ensures that $I$ is actually empty and the two choices must be  $(\xi_1, \dots, \xi_d)$ and $(-\xi_1, \dots, -\xi_d)$, which clearly give rise to convex polytopes.
\end{proof}

Thus we see that the equivariant spectrum of a toric orbifold with no parallel facets and no subpolytopes determines two collections of facet normals and corresponding facet volumes. The question we are asking becomes a purely combinatorial one. 
\begin{question}
Do the normal directions to the facets of a rational simple polytope and the corresponding facet volumes determine the rational simple polytope uniquely?
\end{question}
For convex polytopes, this question and its answer are known as the Minkowski problem.  Daniel Klain \cite{klain} recently gave an elegant solution to this problem.
\begin{thm}\cite[Thm. 2]{klain} \label{klain}
Given a list 
$\{(u_i,\nu_i), u_i  \in \bbR^n, v_i\in \bbR^+,\,\,\,i=1,\dots, d\}$
where the $u_i$ are unit vectors that span $\bbR^n$, there exists a convex polytope $P$ with facet normals $u_1, \dots, u_d$ and corresponding facet volumes $\nu_1, \dots, \nu_d$ if and only if
\[
\sum_{i=1}^d \nu_i u_i = 0.
\]
Moreover, this polytope is unique up to translation.
\end{thm}

Klain proves uniqueness of the solution to the Minkowski problem using a clever inductive argument involving the Minkowski and Brunn-Minkowski inequalities.  These inequalities are generalized isoperimetric inequalities for so-called \emph{mixed volumes}, which encode the relationship between a compact convex set and its orthogonal projections onto subspaces.  For the base case in his induction, Klain gives the solution to the Minkowski problem in dimension $2$ and derives the inequalities as consequences.  He then assumes that the Minkowski and Brunn-Minkowski inequalities hold in dimension $n-1$ and uses them to prove uniqueness in dimension $n$.  The Minkowski inequality in dimension $n$ follows, and this in turn implies the Brunn-Minkowski inequality in dimension $n$.  See \cite{klain} for more details.

Together with Proposition \ref{spectral_data} and Lemma \ref{two_choices}, Theorem \ref{klain} implies the following preliminary result.
\begin{cor}
Let $(\orb,\omega)$ be a toric orbifold with a toric K\"ahler metric such that the moment polytope of $\orb$ has no parallel facets and no subpolytopes.  Then the equivariant spectrum of $\orb$ determines the moment polytope of $\orb$ up to translation and $2$ choices, and hence determines $\orb$ up to symplectomorphism.
\end{cor}

In order to prove Theorem \ref{main_theorem}, we only need to show that rational simple polytopes with no parallel facets and no subpolytopes are generic among all rational simple polytopes. 
\begin{lemma}
Close to any rational simple polytope in $\bbR^n$, there is a rational simple polytope that has no parallel facets and has no subpolytopes.
\end{lemma}

\begin{proof}
Let  $P=\{x \in \bbR^n: x\ip u_i\geq c_i, i=1,\dots, d\}$ be a rational simple polytope and let $\nu_1,\dots, \nu_d$ be the corresponding facet volumes. 

It is easy to perturb our polytope to a rational simple polytope without parallel facets, and this is precisely the advantage of toric orbifolds over toric manifolds. Suppose that $P$ has two parallel facets, say those with labels $1$ and $i$. Choose $\tilde u_1 \in \bbQ^n$ very close to $u_1$. Let 
\[
\tilde{P}=\{x \in \bbR^n: x\ip u_i\geq c_i, i=2,\dots d, \,\,\,x\ip \tilde{u}_1\geq c_1\}.
\]
It is easy to see that $\tilde{P}$ is simple and rational: one can multiply $\tilde{u}_1$ by an appropriate integer multiple so that $\tilde{u}_1 \in \bbZ^n$, and if $\tilde{u}_1$ is sufficiently close to $u_1$, then $\{\tilde{u}_1,u_2,\dots,u_n\}$ will remain a basis of $\bbQ^n$. 

Now we show that one can also perturb our polytope to a rational simple polytope without subpolytopes. If $P$ has a subpolytope then there is a proper subset $I$ of $\{1,\dots,d\}$ such that 
\begin{displaymath}
\sum_I \nu_i u_i=0.
\end{displaymath}
Choose $j \in \{1, \dots, d\} \setminus I$ and perturb $P$ by moving the facet $F_j$ perpendicular to $u_j$ along $u_j$ by distance $\epsilon$, thus changing $c_i$. For $\epsilon$ sufficiently small, we will not introduce new subpolytopes. The only $\nu_i$ 's that are changed by this perturbation are the ones corresponding to facets which intersect $F_j$, and one obtains a one parameter family of polytopes. If the sum
\begin{displaymath}
\sum_I \nu_i(t) u_i
\end{displaymath}
is nonzero for some small value of the parameter $t$, then we have the desired perturbation. If not, then taking derivatives with respect to $t$ yields
\begin{displaymath}
\sum_{I'} \nu_i'(t) u_i=0,
\end{displaymath}
where $I'$ is the subset of $i\in I$ such that $F_i\cap F_j \neq \emptyset$. It is not hard to see that 
\begin{displaymath}
 \nu_i'(0)=\text{Vol}(F_i\cap F_j),
\end{displaymath}
so that
\begin{displaymath}
\sum_{I'} \text{Vol}(F_i\cap F_j)u_i=0.
\end{displaymath}
Choose $k \in I'$ such that $k \neq j$.  By slightly perturbing $u_k$ in $\bbQ^n$ without changing the volume $F_k\cap F_j$ (i.e., such that $u_j^\perp\cap u_k^\perp$ is unchanged), the preceding sum equals a nonzero vector and one obtains the desired rational simple polytope without subpolytopes.
\end{proof}
 
The two possibilities in Theorem \ref{main_theorem} corresponding to the polytopes $P$ and $-P$ may not be equivariantly symplectomorphic, and thus the equivariant spectrum does not determine a single pair (toric orbifold, torus action).  However, the underlying toric manifolds \emph{are} symplectomorphic (but not equivariantly so). 
\begin{prop}\label{unique_symplectic}
Let $(\orb,\omega)$ be a generic toric orbifold endowed with a toric K\"ahler metric. Then the equivariant spectrum of $\orb$ determines the symplectomorphism type of $\orb$. 
\end{prop}
\begin{proof}
Since we have proved Theorem \ref{main_theorem}, we only need to show that the orbifolds $\orb_P$ and $\orb_{-P}$ are symplectomorphic. Let $P$ be given by
\begin{equation} \label{eqn:P}
 P=\{x\in \bbR^n: x\ip m_i u_i\geq c_i,\,\, i=1,\ldots,d\}
\end{equation}
where the $u_i$ are the primitive inward-pointing normals to the facets of $P$ with corresponding weights $m_i$, and $c_i \in \bbR$.
Then $-P$ is given by \eqref{eqn:P} with $-u_i$ replacing $u_i$. The explicit construction of the toric symplectic orbifolds associated to such polytopes (e.g., \cite[pp. 8-9]{m3}) shows that these two orbifolds are in fact the same orbifold, with different torus actions. For the sake of completeness we briefly describe the construction. For the orbifold associated with the polytope $P$ consider the exact sequences
\begin{eqnarray*}\label{exactseq}
0 \rightarrow N \rightarrow \bbT^d \xrightarrow{\beta'} \bbT^n \rightarrow 0 \\
0  \rightarrow  \mathfrak{n} \xrightarrow{\iota} \bbR^d \xrightarrow{\beta} \bbR^n \rightarrow 0
\end{eqnarray*}
where $\beta$ takes an element $e_i$ in the canonical basis for $\bbR^d$ to $m_iu_i$, and $\mathfrak{n}$ is the Lie algebra of $N$.  The group $N$ acts symplectically on $\bbC^d$ with moment map
\begin{displaymath}
 \phi(z)=\sum|z_i|^2\iota^*e_i.
\end{displaymath}
The toric orbifold associated to $P$ is $\orb=\phi^{-1}(c)/N$, 
where $c \in \mathfrak{n}$ is determined by the $c_i$ as $c=-\sum c_i\iota^*e_i$. The symplectic structure on $\orb$ comes from the canonical symplectic structure on $\bbC^d$ via symplectic quotient in the usual way.  We see that this construction will yield the same orbifold $\orb$ with the same symplectic form when we replace $\beta$ by the map $-\beta$ defined by $-\beta(e_i) = -m_i u_i$.  Hence $\orb_P$ and $\orb_{-P}$ are symplectomorphic. 
\end{proof}

Note that the $\bbT^n$-actions on $\orb_P$ and $\orb_{-P}$ differ. For example, for $\orb_P$ we have
\begin{displaymath}
 e^{itu_l}\ip [z_1,\ldots, z_d]=[z_1,\ldots,e^{it}z_l,\ldots,z_d],
\end{displaymath}
whereas we have 
\begin{displaymath}
 e^{itu_l}\ip [z_1,\ldots, z_d]=[z_1,\ldots,e^{-it}z_l,\ldots,z_d]
\end{displaymath}
for $\orb_{-P}$.


\section{Constant scalar curvature metrics are audible}\label{sec:csc}

The goal of this section is to prove Theorem \ref{thm:csc}.
In this theorem, ``generic'' has the same meaning as in Theorem \ref{main_theorem}.   
In particular, the equivariant spectrum of $(\orb,g)$ must determine the moment polytope of $\orb$ (up to translation and two choices) as well as its labels. Theorem \ref{thm:csc} also holds in the setting of smooth toric $2n$-manifolds, provided the corresponding Delzant polytope is audible (again, up to translation and two choices). From \cite{dgs} we know that this occurs for $n=2$ if the manifold is generic and has at most $3$ pairs of parallel sides. 

Throughout this section we will make use of Chern classes for orbifolds. For the relevant background, see \cite[\S2]{gRR}. As for manifolds, Chern classes for orbifolds are diffeomorphism invariants.

We start with a few preliminary results. 
\begin{lemma}\label{lem:charclasses}
The quantities $\int_\orb c_2\wedge \omega^{n-2}$ and  $\int_\orb c_1^2\wedge \omega^{n-2}$ are determined by the equivariant spectrum of $\orb$, as is the cohomology class of $\omega$.
\end{lemma}
Henceforth we will write $\int_\orb c_2$ and $\int_\orb c_1^2$ for $\int_\orb c_2\wedge \omega^{n-2}$ and $\int_\orb c_1^2\wedge \omega^{n-2}$, respectively.
\begin{proof}
We know from \S \ref{sec:mainthm} that the equivariant spectrum determines a set of two rational simple polytopes $P$ and $-P$ that can arise as moment map images of $\orb$. It also determines a unique set of labels $L$ associated to the faces of $P$ (or $-P$). The data $(P,L)$ determine a symplectic toric orbifold $\orb_P$ and the data $(-P, L)$ determine another symplectic toric orbifold $\orb_{-P}$. We know that $\orb$ is equivariantly symplectomorphic to either $\orb_P$ or $\orb_{-P}$. Proposition \ref{unique_symplectic} ensures that  $\orb_P$ and $\orb_{-P}$ are symplectomorphic even though they are not equivariantly symplectomorphic. This then implies that the characteristic classes of $\orb$ are determined from its equivariant spectrum, and so is the cohomology class of $\omega$. 
\end{proof}

Since the labels of the moment polytope of $\orb$ are uniquely determined, the isotropy groups of points in $\orb$ are uniquely determined.
\begin{lemma}\label{hear_iso}
The isotropy groups associated to all points in $\orb$ are determined by the equivariant spectrum. 
\end{lemma}
\begin{proof}
A point $p\in \orb$ is associated to certain facets of its moment polytope $P$; let $\mathcal{F}(p)$ denote the set of facets of $P$ containing $\phi(p)$.  For each facet $F_i$ in $\mathcal{F}(p)$, let $u_i$ denote its primitive outward normal and $m_i$ its label.
Define
\[
\Lambda_p=\text{Span}_{\bbZ}\{u_i: F_i\in \mathcal{F}(p)\},\,\,\,\hat{\Lambda}_p=\text{Span}_{\bbZ}\{m_iu_i: F_i\in \mathcal{F}(p)\}.
\]
In \cite{LT}, Lerman and Tolman show that $\text{Iso}(p)$ is $\Lambda_p/\hat{\Lambda}_p$. Thus $P$ and $L$ determine the isotropy groups of all points in $\orb$, and one can check that the above construction gives the same isotropy groups for $-P$. 
\end{proof}

\begin{rem}
Note that the above characterization of isotropy groups shows that all elements in the pre-image of an open facet in the moment polytope have the same isotropy group. In fact we see that the orbifold stratification is as follows. The highest-dimensional stratum $\mathcal{S}_0$ is the pre-image of the interior of the polytope. The strata of codimension $1$, $\mathcal{S}_{1}$, are the pre-image of the union of the interior of codimension $1$ facets. In general, the codimension $i$ strata, $\mathcal{S}_{i}$, are the pre-image of the union of the relative interior of the faces of codimension $i$.
\end{rem} 

We now prove Theorem \ref{thm:csc}.
\begin{proof}
Let us first apply Theorem \ref{main_asymptotic} to $\orb$ with $f$ equal to the identity. This is the case that is treated in \cite{dggw}. The heat trace is asymptotic to
\[
\sum_{S \in \mathcal{S}(\orb)} \frac{1}{|\text{Iso}(S)|} (4\pi t)^{-\frac{\text{dim}(S)}{2}} \int_{ S} \sum_{k=0}^{\infty} b_k(I, S) t^k \text{dvol}_{S}(x).
\]
In this expansion the term in $t^{-n+2}$ is spectrally determined. It is given by
\begin{equation}\label{eqn:id_exp}
(4\pi)^{-n}b_2(\orb)+(4\pi)^{-n+1}\sum_{i=1}^d \frac{\int_{F_i} b_1(F_i)}{m_i}+(4\pi)^{-n+2}\sum_{a\in I} \frac{\int_{F_a} b_0(F_a)}{m_a},
\end{equation}
where $F_i$ denotes the pre-image via the moment map of a face and $m_i$ is the corresponding label; $I$ denotes some set that indexes codimension $2$ faces; $F_a$ denotes the pre-image via the moment map of a codimension $2$ face; and $m_a=|\is(F_a)|$. 
Note that we use $F_k$ to denote both the $k$th face and its pre-image under the moment map. Although this is an abuse of notation, it should not cause confusion.  We show that each of the terms in \eqref{eqn:id_exp} is determined by the equivariant spectrum. 

We begin with the last term.  We have that $\int_{F_a}b_0(F_a)=\text{Vol}(F_a)$ is determined by $P$ (equivalently, by $-P$), hence it is determined by the equivariant spectrum.  By Lemma \ref{hear_iso} we also hear $m_a$, implying that the last term in \eqref{eqn:id_exp}  is determined by the equivariant spectrum. 

Next we consider the middle term in \eqref{eqn:id_exp}.  
Let $F$ be a face and $u$ the corresponding normal. Consider the isometry $\psi_u$ defined by
\[
\psi_u(p)=e^{iu}.p \ ;
\]
it follows from Theorem \ref{main_asymptotic} that the asymptotic behavior of the $\psi_u$-invariant heat trace is given by
\[
\sum_{S \in \mathcal{S}(\orb)} \frac{1}{|\text{Iso}(S)|} (4\pi t)^{-\frac{\text{dim}(\text{Fix} \psi_u \cap S)}{2}} \int_{\text{Fix} \psi_u \cap S} \sum_{k=0}^{\infty} b_k(f,S) t^k \text{dvol}_{\text{Fix} \psi _u \cap S}(x).
\]
Now $\opF \psi_u=F$, which has dimension $2(n-1)$ and thus does not intersect $\mathcal{S}_0$. The coefficient of $t^{-n+2}$ in this expansion is given by
\[
(4\pi)^{-n+1}\frac{\int_F b_1(F)}{m}+(4\pi)^{-n+2}\sum_{a\in I} \frac{\int_{F_a} b_0(F_a)}{m_a},
\]
where $m$ is the label corresponding to $F$; $I$ denotes some set that indexes the codimension $2$ faces that intersect $F$; $F_a$ denotes the pre-image via the moment map of a codimension $2$ face;  and $m_a=|\is(F_a)|$. 
The same argument as above shows that each summand in the second term of this expansion is determined by the equivariant spectrum.  The index set $I$ is also determined by $P$, so that the sum over $I$ is spectrally determined. Thus we hear $\frac{\int_F b_1(F)}{m}$ for each face, and hence the middle term in \eqref{eqn:id_exp} is determined by the equivariant spectrum. 

Since the middle and last terms in \eqref{eqn:id_exp} are spectrally determined, so is $b_2(\orb)$.  
It follows from \eqref{eqn:b2} that
\begin{equation}\label{b_2}
360 b_2(\orb)= \int _\orb b_2(x) \text{dvol}_{\orb}(x) =\int_\orb  (2|R|^2-2|\rho|^2+5s^2) \text{dvol}_{\orb}(x)
\end{equation}
where $R$ denotes the full curvature tensor, $\rho$ denotes the Ricci tensor, and $s$ denotes the scalar curvature.  Our next goal is to show how the right side of \eqref{b_2} can be expressed as a linear combination of $\int_\orb s^2$, $\int_\orb c_1^2$, and $\int_\orb c_2$, where the coefficient in  $\int_\orb s^2$ is nonzero.

We first show that $\int_\orb |\rho|^2$ can be written as a linear combination of $\int_\orb s^2$ and $\int_\orb c_1^2$. Recall that the complex dimension of $\orb$ is $n$. Let $\rho_0$ be the primitive part of the Ricci curvature so that
\[
\rho=\frac{\text{Tr}(\rho)}{n} \omega+\rho_0.
\]
By definition, $\text{Tr}(\rho)=s$ and the decomposition above is orthogonal, implying that 
\[
|\rho|^2=\frac{s^2}{n^2}|\omega|^2+|\rho_0|^2.
\] 
Since $|\omega|^2=n$ this becomes 
\[
|\rho|^2=\frac{s^2}{n}+|\rho_0|^2.
\]
Using the Apte formula (see \cite[p. 80]{besse}) we have
\begin{equation}\label{eqn:c1squared}
\frac{4\pi^2}{(n-2)!}\int_\orb c_1^2=\frac{n-1}{4n}\int_\orb s^2 -\int _\orb |\rho_ 0|^2,
\end{equation}
and thus 
\begin{align}\nonumber
\int_\orb |\rho|^2=&\frac{1}{n}\int_\orb s^2+\frac{n-1}{4n}\int_\orb s^2-\frac{4\pi^2}{(n-2)!}\int_\orb c_1^2\\ \nonumber
                      =& \frac{n+3}{4n} \int_\orb s^2-\frac{4\pi^2}{(n-2)!}\int_\orb c_1^2\ . \label{ricci}
\end{align}

Next we write $\int_\orb |R|^2$ as a linear combination of $\int_\orb s^2$, $\int_\orb c_1^2$, and $\int_\orb c_2$. One can view $R$ as an endomorphism of $\Omega^2\orb$ and decompose it as
\[
 R=U+Z+W,
\]
where the tensor $U$ is determined by the scalar curvature, $Z$ is related to the trace-free part of the Ricci tensor, and $W$ is the usual Weyl tensor.
This decomposition is also orthogonal so that $|R|^2=|U|^2+|Z|^2+|W|^2$. We have the relations 
\[
|W|^2=|B_0|^2+\frac{3(n-1)}{n+1}|U|^2+\frac{n-2}{n}|Z|^2, \, |U|^2=\frac{s^2}{4n(2n-1)},\, |Z|^2=\frac{|\rho_0|^2}{n-1},
\]
where $B_0$ is the trace-free part of a tensor that arises in a different decomposition of $R$ (see \cite[p. 77]{besse}).  Thus
\begin{eqnarray}
|R|^2& = & |B_0|^2+\frac{2(2n-1)}{n+1}|U|^2+\frac{2(n-1)}{n}|Z|^2 \nonumber \\
         & = & |B_0|^2+\frac{1}{2n(n+1)}s^2+\frac{2}{n}|\rho_0|^2. \label{eqn:Rsquared}
\end{eqnarray}
Using the Apte formula again gives
\[
\frac{8\pi^2}{(n-2)!}\int_\orb c_2=\frac{n-1}{4(n+1)}\int_\orb s^2 -\frac{2(n-1)}{n}\int _\orb |\rho_ 0|^2+\int_\orb |B_0|^2\ ;
\]
we substitute this expression and that from \eqref{eqn:c1squared} into \eqref{eqn:Rsquared} to get
\begin{align}\nonumber
\int_\orb |R|^2 = & \frac{8\pi^2}{(n-2)!}\int_\orb c_2-\frac{n-1}{4(n+1)}\int_\orb s^2 +\frac{2(n-1)}{n}\int _\orb |\rho_ 0|^2 +\frac{1}{2n(n+1)}s^2+\frac{2}{n} \int _\orb |\rho_ 0|^2\\ \nonumber
              =& \frac{8\pi^2}{(n-2)!}\int_\orb c_2+\frac{2+n-n^2}{4n(n+1)}\int_\orb s^2 +2\int _\orb |\rho_ 0|^2\\ \nonumber
              =& \frac{8\pi^2}{(n-2)!}\int_\orb c_2+\frac{2+n-n^2}{4n(n+1)}\int_\orb s^2 -\frac{8\pi^2}{(n-2)!}\int_\orb c_1^2+\frac{n-1}{2n}\int_\orb s^2 \\ \nonumber
              =&\frac{8\pi^2}{(n-2)!}\int_\orb (c_2-c_1^2)+\frac{1}{4}\int_\orb s^2 \ .
\end{align}

We can replace $\int_\orb |\rho|^2$ and $\int_\orb |R|^2$ in (\ref{b_2}) by the equivalent expressions we have found to get
\[
360 b_2(\orb)=\frac{16\pi^2}{(n-2)!}\int_\orb (c_2-c_1^2)+\frac{1}{2}\int_\orb s^2-\frac{n+3}{2n} \int_\orb s^2+\frac{8\pi^2}{(n-2)!}\int_\orb c_1^2+5\int_\orb s^2,
\]
which simplifies to
\[
360 b_2(\orb)=\frac{8\pi^2}{(n-2)!}\int_\orb (2c_2-c_1^2)+\frac{10n-3}{2n}\int_\orb s^2 \ .
\]
We saw above that $b_2(\orb)$ is determined by the equivariant spectrum, and so are $\int_\orb c_1^2$ and $\int_\orb c_2$ by Lemma \ref{lem:charclasses}; hence $\int_\orb s^2$ is spectrally determined. 

We conclude our proof of Theorem \ref{thm:csc} by giving a characterization of constant scalar curvature that is amenable to our spectral setting.  Consider the integral
\[
 \mathcal{C}(g)=\int_\orb (s-\bar{s})^2\frac{\omega^n}{n!}, 
\]
where $\bar{s}$ denotes the average of the scalar curvature, i.e., $\bar{s}=\frac{\int_\orb s}{\text{Vol}(\orb)}$.
It is known that $\bar{s}$ is determined by the symplectic topology of $\orb$ (see \cite{Don_Calabi}): we have 
\[
 \frac{s\omega^n}{n!}=\frac{2\pi c_1\wedge \omega^{n-1}}{(n-1)!}
\]
so that we may write $\bar{s}$ as
\[
 \bar{s}=\frac{2\pi c_1\wedge [\omega]^{n-1}}{\text{Vol}(\orb)(n-1)!}.
\]
The metric $g$ has constant scalar curvature if and only if $\mathcal{C}(g)$ is zero. One can write $\mathcal{C}(g)$ as
\begin{eqnarray*}
 \mathcal{C}(g)&=&\int_\orb (s^2-2\bar{s}s+\bar{s}^2)\frac{\omega^n}{n!} \\ 
                              &=& \int_\orb s^2\frac{\omega^n}{n!}-2\bar{s}\int_\orb s\frac{\omega^n}{n!}+\bar{s}^2 \text{Vol}(\orb)\\
                              &=&\int_\orb s^2\frac{\omega^n}{n!}-2\bar{s}^2\text{Vol}(\orb)+\bar{s}^2 \text{Vol}(\orb)\\
                               &=&  \int_\orb s^2\frac{\omega^n}{n!}-\bar{s}^2\text{Vol}(\orb) 
\end{eqnarray*}
so that $g$ has constant scalar curvature exactly when 
\[
 \int_\orb s^2\frac{\omega^n}{n!}=\frac{1}{\text{Vol}(\orb)}\left(\frac{2\pi c_1\wedge [\omega]^{n-1}}{(n-1)!}\right)^2\ . 
\]
By Lemma \ref{lem:charclasses}, the expression on the right is determined by the equivariant spectrum.
Hence to determine if $g$ has constant scalar curvature we ``hear'' $\int_\orb s^2$, we ``hear'' $\frac{1}{\text{Vol}(\orb)}\left(\frac{2\pi c_1\wedge [\omega]^{n-1}}{(n-1)!}\right)^2$, and we compare the two quantities. 
\end{proof}

\bibliographystyle{plain}
\bibliography{inv_spec}

\end{document}